\documentclass{amsart}
\usepackage{amsfonts,amsmath,amstext,amsthm,amssymb}
\usepackage{hyperref}
\newtheorem{thm}{Theorem}

\newtheorem{lem}[thm]{Lemma}
\theoremstyle{definition}

\newcommand{\sign }{\rm sign\,}
\newcommand{\supp}{\rm supp\,}
\newcommand{\pv}{\rm p.v. }
\newcommand{\ZR}{\ensuremath{\mathbb R}}
\newcommand{\ZZ}{\ensuremath{\mathbb Z}}

\newcommand{\ZI}{\ensuremath{\mathbb I}}

\begin{document}

\title[On unboundedness of maximal operators]
{On unboundedness of maximal operators
for directional Hilbert transforms}%
\author{G. A. Karagulyan}%
\address{Institute of Mathematics Armenian National Academy of
Sciences Marshal Baghramian ave. 24b, Yerevan, 375019, ARMENIA }

 \curraddr{Yerevan State University, Depart. of Applied Mathematics}%
\email{karagul@instmath.sci.am}%


\date{}
\begin{abstract}
We show that for any infinite set of unit vectors $U$ in $\ZR^2$
the maximal operator defined by
\begin{equation*}
H_Uf(x)=\sup_{u\in U}\bigg|\pv\int_{-\infty }^\infty
\frac{f(x-tu)}{t}dt\bigg|,\quad x\in \ZR^2,
\end{equation*}
is not bounded in $L^2(\ZR^2)$.
\end{abstract}
\maketitle
\begin{section}{Introduction}

For a rapidly decreasing function $f$ and a unit vector $u=(\cos
\theta, \sin\theta ), \theta\in [0,2\pi ],$ we define
\begin{equation*}
H_uf(x,y)=\pv\int_{-\infty }^\infty \frac{f(x-t\cos\theta
,y-t\sin\theta )}{t}dt,
\end{equation*}
which is the one dimensional Hilbert transform along the direction
$u$. It is well known this operator can be extended to a bounded
operator from $L^p(\ZR^2)$ to itself when $1<p<\infty $. In this
paper we study operators
\begin{equation}\label{HU}
H_Uf(x,y)=\sup_{u\in U}|H_uf(x,y)|,
\end{equation}
where $U$ is a set of unit vectors $u$ in $\ZR^2$. Analogous
operators for the maximal functions are properly investigated. The
case of lacunary $U$ first are considered in the papers
\cite{Str1}, \cite{CoFe}, \cite{NSW}. A final result is obtained
by A.~Nagel, E.~M.~Stein and S.~Wainger in \cite{NSW}. They proved
the boundedness of the norms of these operators in $L^p$,
$1<p<\infty, $ for a lacunary $U$. Upper bounds of such operators
depending on the cardinality $\# U$ of the set $U$ considered in
the papers
\cite{Barr1},\cite{Barr},\cite{Cor},\cite{Katz1},\cite{Katz2},\cite{Str2}.
And the definitive estimates due to N.~Katz (\cite{Katz1},
\cite{Katz2}). In the papers he obtained a logarithmic order for
the norms of two different maximal operators depending on $\# U$.
Various generalizations of these results are considered in series
of papers (\cite
{Alf},\cite{ASV1},\cite{ASV2},\cite{KaLa},\cite{SS}).

As for the operator (\ref{HU}) there was no any results, but the
bound
\begin{equation*}
\|H_Uf(x)\|_{L^2}\lesssim \log \# U\, \|f\|_{L^2}.
\end{equation*}
This is an immediate consequence of Menshov-Rademacher theorem
(see \cite{KaSt} or \cite{KaSa}), in spite of the fact that in
Katz's theorems subtle range of ideas are used. It was not even
known whether $H_U$ is bounded in $L^2$ or not for an infinite
lacunary set $U$. The main result of the paper is unboundedness of
$H_U$ in $L^2$ for any infinite $U$.
\begin{thm}\label{T1}
For any infinite set of unit vectors $U$ the operator $H_U$ can
not be extended to a bounded operator from $L^2(\ZR^2)$ to
$L^p(\ZR^2)$ for all $1\le p<\infty $.
\end{thm}
This theorem is an immediate consequence of the following
estimate.
\begin{thm}\label{T2}
If $U$ is a finite set and $1\le p<\infty $ then $\|H_U\|_{L^2\to
L^p}\ge c \sqrt {\log \#U}$, where $c>0$ is an absolute constant.
\end{thm}
The author is very grateful to M.~T.~Lacey for a nice hospitality
during two visits at the Georgia Tech in 2003 and 2004, and for
helpful conversations about these problems.
\end{section}
\begin{section}{Proof of Theorems}
Let $f_n(x)$, $n=1,2,\cdots , 2^m-1$ ($f_n\not\equiv 0$) be a
system of functions defined on the square
\begin{equation*}
Q=[-\pi,\pi ]\times [-\pi,\pi ].
\end{equation*}
In some places we shall use for $f_n(x)$ {\it double numbering},
defining by
\begin{equation*}
f_j^{(k)}(x)=f_n(x),\quad n=2^k+j-1,\quad 1\le j\le 2^k, \quad
k=0,1,\cdots ,m-1.
\end{equation*}
 We shall say the sequence $f_n(x)=f_j^{(k)}(x)$ is {\it
tree-system}, if
\begin{equation}\label{2-25}
\supp f_{2j-1}^{(k+1)}\subset \{ x\in Q: f_j^{(k)}(x)> 0\},\quad
\supp f_{2j}^{(k+1)}\subset \{ x\in Q: f_j^{(k)}(x)< 0\}.
\end{equation}
Applying (\ref{2-25}) several times we get
\begin{equation}
\begin{array}{lll}
 &\supp f_i^{(k+r)}\subset \{ x\in Q: f_j^{(k)}(x)>
0\}&\quad\Leftrightarrow \quad i\in \big(2^rj-2^r, 2^rj-2^{r-1}\big],\\
 &\supp f_i^{(k+r)}\subset \{ x\in Q: f_j^{(k)}(x)<
0\}&\quad\Leftrightarrow \quad i\in\big(2^rj-2^{r-1},2^rj\big],
\end{array}
\end{equation}
and then
\begin{eqnarray}
&\supp f_i^{(k+r)}\cap \{ x\in Q: f_j^{(k)}(x)> 0\}=\varnothing
&\quad\Leftrightarrow \quad
i\not\in \big(2^rj-2^r, 2^rj-2^{r-1}\big],\label{2-1}\\
&\supp f_i^{(k+r)}\cap \{ x\in Q: f_j^{(k)}(x)< 0\}=\varnothing
&\quad\Leftrightarrow \quad
i\not\in\big(2^rj-2^{r-1},2^rj\big].\label{2-2}
\end{eqnarray}
The following lemma for the Haar system is proved by
E.~M.~Nikishin and P.~L.~Ul'yanov \cite {NiUl} and we use the same
idea to prove a general one (see also \cite {KaSa}).
\begin{lem}\label{L1}
If $f_n(x)$, $x\in Q$, $n=1,2,\cdots ,2^m-1 $, is tree-system then
there exists a permutation $\sigma $ of the numbers $\{1,2,\cdots
,2^m-1\}$ such that
\begin{equation*}
\sup_{1\le l< 2^m}\big|\sum_{n=1}^{l}f_{\sigma (n)}(x)\big|\ge
\frac{1}{3}\sum_{n=1}^{2^m-1}\big|f_n(x)\big|.
\end{equation*}
\end{lem}
\begin{proof}
We connect with each $f_n(x),\, n=2^k+j-1$, a number
\begin{equation*}
t_n=\frac{2j-1}{2^{k+1}}\in [0,1].
\end{equation*}
 Notice they are not equal for different $n$'s. Define the permutation
 $\sigma $ so that
\begin{equation*}
t_{\sigma (1)}< t_{\sigma (2)}< \cdots < t_{\sigma (2^m-1)}.
\end{equation*}
We shall prove that for any $x\in Q$ there exists a number
$l=l(x)$ with
\begin{gather}
f_{\sigma (n)}(x)\ge 0,\hbox { if } n>l(x),\label{lx1}\\
f_{\sigma (n)}(x)\le 0,\hbox { if } n\le l(x).\label{lx2}
\end{gather}
Defining
\begin{equation*}
l=l(x)=\sup \{ n:\, 1\le n< 2^m,\,f_{\sigma (n)}(x)\le 0\},
\end{equation*}
we shall have (\ref{lx2}) immediately and if $\nu =\sigma
(l(x)+1)$ then
\begin{equation}\label{2-3}
f_\nu (x)>0.
\end{equation}
To prove (\ref{lx1}) it is enough to show that if (\ref{2-3})
holds and $t_n>t_\nu $ then $f_n(x)\ge 0$ . Suppose
\begin{equation*}
\nu=2^k+j-1,\quad n=2^s+i-1.
\end{equation*}
According to the assumption
\begin{equation}\label{2-30}
t_\nu=\frac{2j-1}{2^{k+1}}<\frac{2i-1}{2^{s+1}}=t_n.
\end{equation}
If $s> k$ then $s=k+r\,(r>0)$. From (\ref{2-30}) we get
\begin{equation*}
i>2^rj-2^{r-1}.
\end{equation*}
Therefore by (\ref{2-1}) we obtain
\begin{equation*}
\supp f_n\cap \{f_\nu >0\}=\varnothing .
\end{equation*}
Since $f_\nu (x)>0$ we get $f_n(x)=0$.

\noindent If $s<k$, then $k=s+r\,(r>0)$. Applying (\ref{2-30}) we
get
\begin{equation*}
j\le 2^ri-2^{r-1}.
\end{equation*}
Hence by (\ref{2-2}) we have
\begin{equation*}
\supp f_\nu \cap \{f_n<0\}=\varnothing .
\end{equation*}
By (\ref{2-3}) we conclude $f_n(x)\ge 0$. So (\ref{lx1}) and
(\ref{lx2}) are proved. Using them, we obtain
\begin{equation*}
\max_{1\le l< 2^m}\big|\sum_{n=1}^{l}f_{\sigma (n)}(x)\big|=\max\{
-\sum_{n=1}^{l(x)}f_{\sigma (n)}
(x),\big|\sum_{n=1}^{2^m-1}f_{\sigma (n)}(x)\big|\}.
\end{equation*}
If
\begin{equation*}
-\sum_{n=1}^{l(x)}f_{\sigma (n)} (x)=\sum_{n=1}^{l(x)}|f_{\sigma
(n)}(x)|<\frac{1}{3} \sum_{n=1}^{2^m-1}|f_{\sigma (n)}(x)|
\end{equation*}
then we get
\begin{equation*}
\sum_{n=l(x)+1}^{2^m-1}|f_{\sigma (n)}(x)|>
2\sum_{n=1}^{l(x)}|f_{\sigma (n)}(x)|
\end{equation*}
and therefore
\begin{multline*}
\big|\sum_{n=1}^{2^m-1}f_{\sigma (n)}(x)\big|\ge
|\sum_{n=l(x)+1}^{2^m-1}f_{\sigma (n)}(x)|-
|\sum_{n=1}^{l(x)}f_{\sigma
(n)}(x)|\\
>\frac{1}{2}\sum_{n=l(x)+1}^{2^m-1}|f_{\sigma (n)}(x)|>
\frac{1}{3}\sum_{n=1}^{2^m-1}|f_{\sigma (n)}(x)|.
\end{multline*}
Thus we conclude
\begin{equation*}
\sup_{1\le l< 2^m}\big|\sum_{n=1}^{l}f_{\sigma
(n)}(x)\big|>\frac{1}{3}\sum_{n=1}^{2^m-1}|f_{\sigma
(n)}(x)|=\frac{1}{3}\sum_{n=1}^{2^m-1}|f_n(x)|.
\end{equation*}
\end{proof}
Fix a Schwartz function $\phi (x)$ with
\begin{equation}
\phi (x) > 0, \quad  \int_\ZR\phi (x)dx=1,\quad  \supp\widehat\phi
\subset[-1 ,1].\label{phi}
\end{equation}
We consider operators
\begin{equation}
\Phi_n(f)=\Phi_n(x,y,f)=n^2\int_\ZR
f(x-t,y-s)\phi(nt)\phi(ns)dtds,\quad n=1,2,\cdots .
\end{equation}
Applying (\ref{phi}), for any $f\in L^\infty (\ZR^2)$ we get
\begin{gather}
 \inf_{(x,y)\in \ZR^2}f(x,y)\le \Phi_n(x,y,f) \le \sup_{(x,y)\in
 \ZR^2}f(x,y),\label{Phi}\\
\supp\widehat{\Phi_n(f)}\subset [-n,n]\times [-n,n].\label{2-5}
\end{gather}
 If in addition $f$ is compactly supported, using a standard argument, we
conclude
\begin{equation}\label{2-4}
\|\Phi_n(f)-f\|_{L^2}\to 0,\hbox { as } n\to \infty .
\end{equation}
If
\begin{equation}\label{2-9}
n=2^k+j-1,\quad 1\le j\le 2^k, \quad k=0,1,\cdots ,m-1,
\end{equation}
then we denote
\begin{equation*}
\bar n=2^{k-1}+\bigg[\frac{j+1}{2}\bigg],
\end{equation*}
where $[\cdot ]$ means an integer part of a number. Using this
notation we may write the conditions (\ref{2-25}) by
\begin{equation}\label{2-17}
\supp f_n\subset \{(-1)^{j+1}\cdot f_{\bar n}>0\}.
\end{equation}
We shall consider sectors defined by
\begin{equation*}
\{(x,y)\in \ZR^2:\, x+iy=re^{i\theta },\, r\ge 0,\, \alpha\le
\theta \le \beta \}
\end{equation*}
where $0\le \alpha <\beta \le 2\pi $. Some arguments in the prove
of following lemma are coming from the paper \cite {Kar}.
\begin{lem}\label{L2}
Let $S_n$, $n=1,2,\cdots ,\nu =2^m-1$, be sectors on the plane.
Then there exist functions $f_n\in L^2(\ZR^2)$, $n=1,2,\cdots ,\nu
$, such that
\begin{eqnarray}
&\supp \widehat{f_n}\subset S_n,\label{2-6}\\
&\sum_{j=1}^\nu \|f_j\|_{L^2}^2\le c_1,\label{2-7}\\
&|\{(x,y)\in Q:\max_{1\le n\le \nu }|\sum_{j=1}^nf_{\sigma
(j)}(x,y)|>c_3\sqrt{\log \nu }\}|>c_2,\label{2-8}
\end{eqnarray}
where $\sigma $ is the  permutation from Lemma~\ref{L1} and all
the constants are absolute .
\end{lem}
\begin{proof}
We will assume (\ref{2-9}) everywhere below. For a given
$\varepsilon
>0$ define sets $E_n=E_j^{(k)}\subset Q$ ($E_1=E_1^{(0)}=Q$), $g_n\in
L^\infty (\ZR^2)$ and $p_n,q_n\in \ZZ$ with conditions

    a) $E_n=\big\{(x,y)\in E_{\bar n} : (-1)^{j+1}
    \cos (p_{\bar n}x+q_{\bar n}y)>0\big\}$, $n=2,3,\cdots ,\nu  $,

    b) $0\le g_n\le 1$, $\|g_n-\ZI_{E_n}\|_{L^2}\le \varepsilon $, $n=1,2,\cdots ,\nu  $,

    c) $\supp \widehat{g}_n\subset S_n-(p_n,q_n)$, $n=1,2,\cdots
    ,\nu  $,

    d) $\int_{E_n}|\cos (p_nx+q_ny)|dxdy>\frac{|E_n|}{3}$.

\noindent
 We do it by induction. Take $E_1=E_1^{(0)}=Q$.
According to (\ref{2-4}) there exists $l>0$ with
\begin{equation*}
\|\Phi_l(\ZI_{E_1})-\ZI_{E_1}\|_{L^2}<\varepsilon .
\end{equation*}
Define $g_1=\Phi_l(\ZI_{E_1})$ and then applying (\ref{Phi}) we
get b) for n=1. We note that if $E$ is a measurable set then
\begin{multline}\label{2-14} \int_E|\cos (px+qy)|dxdy\ge
\int_E\cos^2(px+qy)dxdy\\
=\frac{|E|}{2}+ \int_E\frac{\cos(2(px+qy))}{2}dxdy\rightarrow
\frac{|E|}{2}\hbox { as } |p|,|q|\to\infty .
\end{multline}
This observation shows that for sufficiently large $p_1=p$ and
$q_1=q$ we shall have condition d) for $n=1$. On the other hand by
(\ref{2-5}) $\supp \widehat{g}_1$ is bounded.  Thus for an
appropriate $p_1,q_1$ we will have also c) (with $n=1$). Certainly
we can choose $p_1$ and $q_1$ common for both conditions c) and
d). Now we suppose that the conditions a)-d) holds for any $k<n$,
in particular for $\bar n$. We define $E_n$ by the equality in a).
Then we chose positive integer $l$ with
\begin{equation*}
\|\Phi_l(\ZI_{E_n})-\ZI_{E_n}\|_{L^2}<\varepsilon  .
\end{equation*}
and denote $g_n=\Phi_l(\ZI_{E_n})$. Again applying (\ref{2-14})
and using the boundedness of $\supp \widehat{g}_n$ we may chose
integers $p_n,q_n$ satisfying c) and d) together. Using the
condition a), it is easy to check that the sets $E_n$ satisfies
the conditions
\begin{gather*}
E_j^{(k)}\cap E_{j'}^{(k)}=\varnothing,\hbox { if } j\neq j',\\
E_{2j-1}^{(k+1)}\cup E_{2j}^{(k+1)}\subset E_j^{(k)},\quad
|E_j^{(k)}\setminus (E_{2j-1}^{(k+1)}\cup E_{2j}^{(k+1)})|=0.
\end{gather*}
Using this properties we conclude
\begin{equation}\label{2-21}
\sum_{n=1}^\nu \ZI_{E_n}(x,y)=m,\hbox { almost everywhere on } Q.
\end{equation}
 Now we define
\begin{equation}\label{2-16}
f_n(x,y)=\frac{e^{i(p_nx+q_ny)}g_n(x,y)}{\sqrt m}.
\end{equation}
The condition (\ref{2-6}) immediately follows from c), because
\begin{equation*}
\supp \widehat{f}_n=\supp \widehat{g}_n+(p_n,q_n).
\end{equation*}
On the other hand taking a small $\varepsilon $ by b) and
(\ref{2-21}) we obtain
\begin{equation*}
\sum_{n=1}^\nu \int_Q|f_n|^2=\frac{1}{m}\sum_{n=1}^\nu
\int_Q|g_n|^2\le \frac{2\nu
}{m}\varepsilon^2+\frac{2}{m}\sum_{n=1}^\nu \int_Q\ZI_{E_n}\le
c_1,
\end{equation*}
which gives (\ref{2-7}). Now consider functions
\begin{equation}\label{2-13}
\tilde f_n=\Re f_n\cdot \ZI_{E_n}=\frac{\cos(p_nx+q_ny)\cdot
g_n(x,y)\cdot \ZI_{E_n}(x,y)}{\sqrt m}.
\end{equation}
 Applying b) and (\ref{2-16}) we get
\begin{equation}\label{2-10}
\|\tilde f_n-\Re f_n\|_{L^2}^2=\int_{\ZR^2\setminus E_n}|\Re
f_n|^2\le\int_{\ZR^2\setminus E_n}|f_n|^2\\
= \frac{1}{m}\int_{\ZR^2\setminus E_n}|g_n|^2\le \frac{\varepsilon
^2}{m}.
\end{equation}
On the other hand we have
\begin{equation}\label{2-11}
\max_{1\le n\le N}|\sum_{j=1}^nf_{\sigma (j)}|\ge \max_{1\le n\le
N}|\sum_{j=1}^n\Re f_{\sigma (j)}| \ge\max_{1\le n\le
N}|\sum_{j=1}^n\tilde f_{\sigma (j)}|-\sum_{j=1}^\nu |\tilde
f_j-\Re f_j|.
\end{equation}
From (\ref{2-10}) we obtain
\begin{equation*}
\bigg\|\sum_{j=1}^\nu |\tilde f_j-\Re f_j|\bigg\|_{L^2}\le
\frac{\nu \varepsilon }{\sqrt m}.
\end{equation*}
Therefore taking a small $\varepsilon >0$ we can say that
\begin{equation}\label{2-12}
\bigg|\bigg\{(x,y)\in Q:\, \sum_{j=1}^\nu |\tilde f_j-\Re
f_j|>1\bigg\}\bigg|\le \delta,
\end{equation}
for any given $\delta >0$. From (\ref{2-11}) and (\ref{2-12}) we
conclude, that to prove (\ref{2-8}) and so the lemma it is enough
to prove
\begin{equation}\label{2-22}
|\{(x,y)\in Q:\max_{1\le n\le \nu }|\sum_{j=1}^n\tilde f_{\sigma
(j)}(x,y)|>c_3\sqrt{\log \nu  }\}|>c_2
\end{equation}
Let us show that $\tilde f_n$ is a  tree-system, i.e. it satisfies
(\ref{2-17}). Since $g_n> 0$ from (\ref{2-13}) we get that $\tilde
f_n(x,y)$ and $\cos(p_nx+q_ny)$ have same sign in the set $E_n$.
Therefore by a) we obtain
\begin{multline*}
\supp \tilde f_n\subset E_n =\{(x,y)\in E_{\bar n} :\, (-1)^{j+1}
    \cos (p_{\bar n}x+q_{\bar n}y)>0\big\}\\
    =\{(x,y)\in E_{\bar n} :\,
    (-1)^{j+1}\tilde f_{\bar n}(x,y)>0\big\}=\{(x,y)\in E_{\bar n} :\,
    (-1)^{j+1}\tilde f_{\bar n}(x,y)>0\big\}
    \\
    =\{(x,y)\in Q:\,(-1)^{j+1}\tilde f_{\bar n}(x,y)>0\big\}.
\end{multline*}
Hence $\tilde f_n$ is tree-system. So according to Lemma~\ref{L1}
we have
\begin{equation}\label{2-18}
\max_{1\le n\le N}|\sum_{j=1}^n\tilde f_{\sigma
(j)}(x,y)\ge\frac{1}{3}\sum_{j=1}^\nu |\tilde f_j(x,y)|.
\end{equation}
From (\ref{2-16}) and the conditions b) and d) we get
\begin{multline*}
\int_Q|\tilde f_n|=\frac{1}{\sqrt
m}\int_{E_n}|g_n\cos(p_nx+q_ny)|dxdy \ge\frac{1}{\sqrt m}
\int_{E_n}|\cos(p_nx+q_ny)|dxdy
\\
-\frac{1}{\sqrt m} \int_{E_n}|(1-g_n)\cos(p_nx+q_ny)|dxdy\ge
\frac{|E_n|}{3\sqrt m}-\frac{\varepsilon}{\sqrt m}.
\end{multline*}
If we take $\varepsilon >0$ to be small then from (\ref{2-21}) we
obtain
\begin{equation*}
\int_Q\sum_{j=1}^\nu |\tilde f_j|\ge \frac{1}{3\sqrt
m}\sum_{n=1}^\nu|E_n|-\frac{\nu \varepsilon }{\sqrt m}\gtrsim
\sqrt m.
\end{equation*}
Combining this and (\ref{2-18}) we get
\begin{equation}\label{2-19}
\int_Q\max_{1\le n\le N}\big|\sum_{j=1}^n\tilde f_{\sigma
(j)}(x,y)\big|\gtrsim \sqrt m .
\end{equation}
On the other hand by (\ref{2-21}), (\ref{2-16}) and b) for any
$(x,y)\in Q$ we have
\begin{equation}\label{2-20}
\max_{1\le n\le N}|\sum_{j=1}^n\tilde f_{\sigma (j)}(x,y)|\le
\sum_{j=1}^\nu |\tilde f_j(x,y)|\le \frac{1}{\sqrt
m}\sum_{j=1}^\nu \ZI_{E_n}(x,y)\le \sqrt m.
\end{equation}
From (\ref{2-19}) and (\ref{2-20}) follows (\ref{2-22}).
\end{proof}
\begin{proof}[Proof of Theorem~\ref{T2}] For any region $S\subset
\ZR^2$ we denote
\begin{equation*}
T_Sf(x,y)=\int_S e^{i(\xi x+\eta y)}\widehat{f}(\xi ,\eta )d\xi
d\eta.
\end{equation*}
Since the multiplier for the Hilbert transform is $i\cdot\sign x$,
for any direction $u=(\cos \theta, \sin\theta )$ we have
\begin{equation*}
\widehat{H_u}f(x,y)=i\cdot \sign (x\cos \theta +y\sin\theta
)\widehat{f}(x,y).
\end{equation*}
Thus we conclude
\begin{equation}\label{2-27}
H_uf=i(2\cdot T_{\Gamma_u}f-f)
\end{equation}
where
\begin{equation*}
\Gamma_u=\{(x,y)\in \ZR^2:\, x\cos \theta +y\sin\theta \ge 0\}.
\end{equation*}
Denote
\begin{equation*}
T_Uf=\sup_{u\in U}|T_{\Gamma_u}f|.
\end{equation*}
Let $U=\{u_k=(\cos \theta_k , \sin\theta_k ):\, k=1,2,\cdots ,N\}$
be the set of directions from Theorem 2. Without loss of
generality we can assume $\theta_k\in (0,\pi/2)$,
$\theta_1<\theta_2<\cdots <\theta_N$ and $N=2^m$. According to
(\ref{2-27}), to prove Theorem 2 it is enough to prove that
\begin{equation}\label{2-28}
\|T_Uf\|_{L^1}\gtrsim \sqrt {\log N}\|f\|_{L^2}
\end{equation}
for some function $f\in L^2(\ZR^2)$. We denote by
$S_k,\,k=1,2,\cdots ,\nu =2^m-1$, the sectors obtained by the
vectors $(u_k)^\perp =(\cos \theta_k , -\sin\theta_k )$ and
$(u_{k+1})^\perp =(\cos \theta_{k+1} , -\sin\theta_{k+1} )$, i.e.
\begin{equation}\label{2-24}
S_k=\{(x,y)\in \ZR^2:\,x\ge 0, x\cos \theta_k +y\sin\theta_k \ge
0,\, x\cos \theta_{k+1} +y\sin\theta_{k+1} \le 0\}.
\end{equation}
Hence if we suppose
\begin{equation*}
\supp\hat f\subset \bigcup_{k=1}^{\nu }S_k
\end{equation*}
then we can write
\begin{equation}\label{2-23}
T_Uf(x,y)=\sup_{1\le l\le
\nu}\bigg|\sum_{k=1}^lT_{S_k}f(x,y)\bigg|.
\end{equation}
We define functions $f_n$ satisfying the conditions of the
Lemma~\ref{L2} corresponding to the sectors $S'_n=S_{\sigma^{-1}
(n)}$ in (\ref{2-24}). Denote
\begin{equation*}
f=\sum_{k=1}^\nu f_k.
\end{equation*}
Since $S_n$ are mutually disjoint the functions $f_n$ are
orthogonal. Thus by (\ref{2-7}) we get $\|f\|_{L^2}\le c_1$. From
(\ref{2-6}) we have $\supp \widehat {f}_n\subset
S_{\sigma^{-1}(n)}$, $n=1,2,\cdots ,\nu $, i.e. $\supp
\widehat{f}_{\sigma (n)}\subset S_n$ and therefore
\begin{equation*}
T_{S_n}f(x)=f_n(x).
\end{equation*}
According to (\ref{2-23}) we obtain
\begin{equation*}
T_Uf(x,y)=\max_{1\le l\le N}|\sum_{j=1}^lf_{\sigma (j)}(x,y)|.
\end{equation*}
Using (\ref{2-7}) and (\ref{2-8}), we get
\begin{equation*}
|\{(x,y)\in Q:T_Uf(x,y)>c_3\sqrt{\log \nu }\}|>c_2
\end{equation*}
And therefore $\|T_Uf\|_{L^p}\gtrsim \sqrt{\log N }\|f\|_{L^2}$.
\end{proof}
\end{section}

\bibliographystyle{amsplain}

\end{document}